\documentclass[reqno,12pt]{amsart}
\usepackage{epsf,epsfig}
\usepackage{graphics,color}
\usepackage{amsmath}
\usepackage{amssymb}
\usepackage{cite}
\usepackage{verbatim}
\usepackage{float}
\usepackage{graphicx}
\usepackage{amsthm}
\usepackage{textcomp}
\usepackage{subfig}
\usepackage{pgf}
\usepackage{ragged2e}

\newtheorem{theorem}{Theorem}[section]

\newtheorem{corollary}[theorem]{Corollary}

\theoremstyle{definition}

\numberwithin{equation}{section}
\renewcommand{\d}{\mathrm{d}}

\newcommand{\epsi}{\varepsilon}
\newcommand{\Rz}{{\mathbb R}}
\newcommand{\Nz}{\mathbb{N}}

\newcommand{\disp}{\displaystyle}
\newcommand{\haz}{\widehat}
\newcommand{\ove}{\overline}

\newcommand{\weak}{\rightharpoonup}
\newcommand{\wstarto}{\stackrel{\ast}{\rightharpoonup}}
\DeclareMathOperator*{\argmin}{arg\,min}
\newcommand{\UUU}{\color{black}}
\newcommand{\EEE}{\color{black}}

\newcommand{\bC}{{\mathbb C}}
\newcommand{\dx}{{\rm d} x}

\begin{document}

\title[Filtered/phase-field method in elasticity]{Analysis of a
  combined Filtered/phase-field approach to topology optimization in elasticity}

\author[F. Auricchio]{Ferdinando Auricchio}
\address[Ferdinando Auricchio]{Department of Civil Engineering and Architecture, University of Pavia, via Ferrata 3, I-27100 Pavia, Italy, \& Istituto di
  Matematica Applicata e Tecnologie Informatiche {\it E. Magenes}, via
  Ferrata 1, I-27100 Pavia, Italy}
\email{auricchio@unipv.it}
\urladdr{http://www-2.unipv.it/auricchio/}

\author[M. Marino]{Michele Marino}
\address[Michele Marino]{Department of Civil Engineering and Computer Science, University of Rome Tor Vergata, Via del Politecnico 1, I-00133 Roma, Italy} 
\email{m.marino@ing.uniroma2.it}

\author[I. Mazari]{Idriss Mazari}
\address[Idriss Mazari]{Universit\'e Paris-Dauphine,
Place du Mar\'echal De Lattre De Tassigny,
F-75775 Paris CEDEX 16, France} 
\email{mazari@ceremade.dauphine.fr}
\urladdr{https://www.ceremade.dauphine.fr/$\sim$mazari/}

\author[U. Stefanelli]{Ulisse Stefanelli} 
\address[Ulisse Stefanelli]{Faculty of Mathematics, University of
  Vienna, Oskar-Morgenstern-Platz 1, A-1090 Vienna, Austria,
Vienna Research Platform on Accelerating
  Photoreaction Discovery, University of Vienna, W\"ahringerstra\ss e 17, 1090 Wien, Austria,
 \& Istituto di
  Matematica Applicata e Tecnologie Informatiche {\it E. Magenes}, via
  Ferrata 1, I-27100 Pavia, Italy
}
\email{ulisse.stefanelli@univie.ac.at}
\urladdr{http://www.mat.univie.ac.at/$\sim$stefanelli/}

 \subjclass[2010]{Topology optimization, elasticity, filter, phase
   field, existence, $\Gamma$-convergence, space-discretization, Lagrangian formulation}

 \begin{abstract}
We advance a combined filtered/phase-field approach to
topology optimization in the setting of linearized
elasticity. Existence of minimizers is proved and rigorous parameter
asymptotics are discussed by means of
variational convergence techniques. Moreover, we investigate an abstract space
discretization in the spirit of conformal finite elements. \UUU
Eventually, stationarity  is equivalently reformulated in terms of a 
Lagrangian. \EEE
 \end{abstract}

\maketitle


\section{Introduction}

Topology optimization is concerned with the determination of optimal
shapes with respect to a given target. In the elastic setting, it often
consists in identifying the portion of a given design domain $\Omega \subset \Rz^3$ (open, smooth,
connected) to be occupied by an elastic solid, so that the
compliance corresponding to its equilibrium state is minimal. The order parameter $\phi:\Omega \to [0,1]$ describes
the presence of material in the domain. In particular, the set $\{\phi=1\}$ is
the solid to be identified, whereas $\{\phi=0\}$ is interpreted as a very
compliant {\it Ersatz} material, still assumed to be elastic. 

The
elastic response is modeled by the continuously differentiable {\it
  elasticity tensor} $\phi \mapsto
\bC(\phi)$, taking values in the symmetric, isotropic 4-tensors, with
$c_0 {\mathbb I}\leq \bC(\cdot)\leq {\mathbb I}/c_0$ for some $c_0>0$,
where ${\mathbb I}$ is the identity 4-tensor and  
${\mathbb A}\leq {\mathbb B}$ if and only if $ {\mathbb B}- {\mathbb
  A}$ is positive semidefinite. \UUU A classical choice for $\bC$ would
be $\bC(\phi) =\bC_0 + \phi^q (\bC_1 - \bC_0) $ with $\bC_1\geq
\bC_0\geq c_0 {\mathbb I} $ and $q\geq 1$. This would correspond to
associate the tensors 
$\bC_0$ and $\bC_1$ to the phases $\phi=0$ and $\phi=1$,
respectively. \EEE
The equilibrium of the body is described by the system
\begin{align}
  \nabla \cdot \bC(\phi) \epsi(u) + \phi f =0 \quad&\text{in} \ \Omega,\label{eq:a}\\
  \bC(\phi) \epsi(u)n =g \quad&\text{in} \ \Gamma_{\rm N}, \\
u=0\quad&\text{in} \ \Gamma_{\rm D}.\label{eq:3}
\end{align}
Here, $\Gamma_{\rm
  D}$ and $\Gamma_{\rm N}$ are two distinct portions (open in the
topology of $\partial \Omega$, such that $\ove\Gamma_{\rm D}\cup \ove
\Gamma_{\rm N} = \partial \Omega$) of the boundary
$\partial \Omega$ where the body is clamped and a traction is exerted,
respectively. In particular, $g:
\Gamma_{\rm N} \to \Rz^3$ is a surface traction density, and $n$ is
the outward unit normal to $\partial \Omega$, while
$f:\Omega \to \Rz^3$ is a force density per unit $\phi$. We assume that $f\in
L^2(\Omega;\Rz^3)$ and $g \in L^2(\Gamma_{\rm N};\Rz^3)$.

Our goal is to minimize the compliance
$$C(\phi,u):= \int_\Omega \phi f \cdot u  \, \dx+\int_{\Gamma_N} g \cdot u \,
      \d \Gamma,$$
with respect to the order parameter $\phi$. Here, $\phi$ is assumed to
take values in $ [0,1]$ almost everywhere with $\int_{\Omega}  \phi(x)\,
      \d x = v_0$, where $0<v_0 <|\Omega|$ is a specified volume, and $u$ is the unique
solution of the equilibrium system \eqref{eq:a}-\eqref{eq:3}, given
$\phi$.  
We can put this {\it topology optimization} problem in variational
terms as the {\it bilevel minimization} problem
\begin{equation}
  \label{eq:murat}
  \min_{\phi \in \Phi} \Big\{C(\phi,u): \ u = \argmin_{v \in U} \big(E(\phi,v)-C(\phi,v)\big)\Big\}.
\end{equation}
Here, $E$ indicates the elastic energy 
\begin{align*}
  E(\phi, u) &=\frac12 \int_\Omega \bC(\phi) \epsi(u):\epsi(u) \, \dx
\end{align*} 
and we have used the notation
\begin{align*}
  \Phi&=\left\{\phi \in L^\infty(\Omega) \ : \ 0\leq \phi \leq
      1 \ \text{a.e.}, \ \int_{\Omega}  \phi(x)\,
        \d x = v_0\right\},\\
  U&:= \{u \in H^1(\Omega;\Rz^3) \ : \ u =0 \ \text{on} \ \Gamma_{\rm
   D}\}
\end{align*}  
for the state spaces.

The topology optimization problem \eqref{eq:murat} cannot be expected to be solvable, since
minimizing sequences may develop fine-scaled 
oscillations. Some examples in this direction are already in 
\cite{Murat}, even in the simpler purely elliptic case. This lack of compactness may be tamed by weakening the
solution concept, namely, by dropping the functional dependence $\phi\mapsto \bC(\phi)$
and considering the order parameter $\phi$ and the elastic strain
$\bC$ as independent variables. This corresponds to the so-called {\it
  homogenization method}, see the classical
monograph \cite{Allaire} for a comprehensive discussion
on its theory and application.

If one is interested in retaining the functional dependence
$\phi\mapsto \bC(\phi)$, problem \eqref{eq:murat} calls for a
regularization. Two prominent possibilities in this direction are the
{\it Filtered} (F) method \cite{Borrvall,Bourdin} and the
{\it Phase-Field} (PF) approach \cite{Bourdin2,Burger}. In the former, the dependence of the
elastic tensor on the material is usually combined with the action of
a {\it filter}, which essentially \UUU amounts \EEE to a regularization. A very
effective and adapted approach within the family of filtered methods
is the so called {\it SIMP} method \cite{BS,BS2}, which is
nowadays often exploited due to its performance and simplicity. 
In the PF approach,  variations of $\phi$ are additionally penalized.

\begin{figure}[h] 
  \centering
\includegraphics[width=\columnwidth]{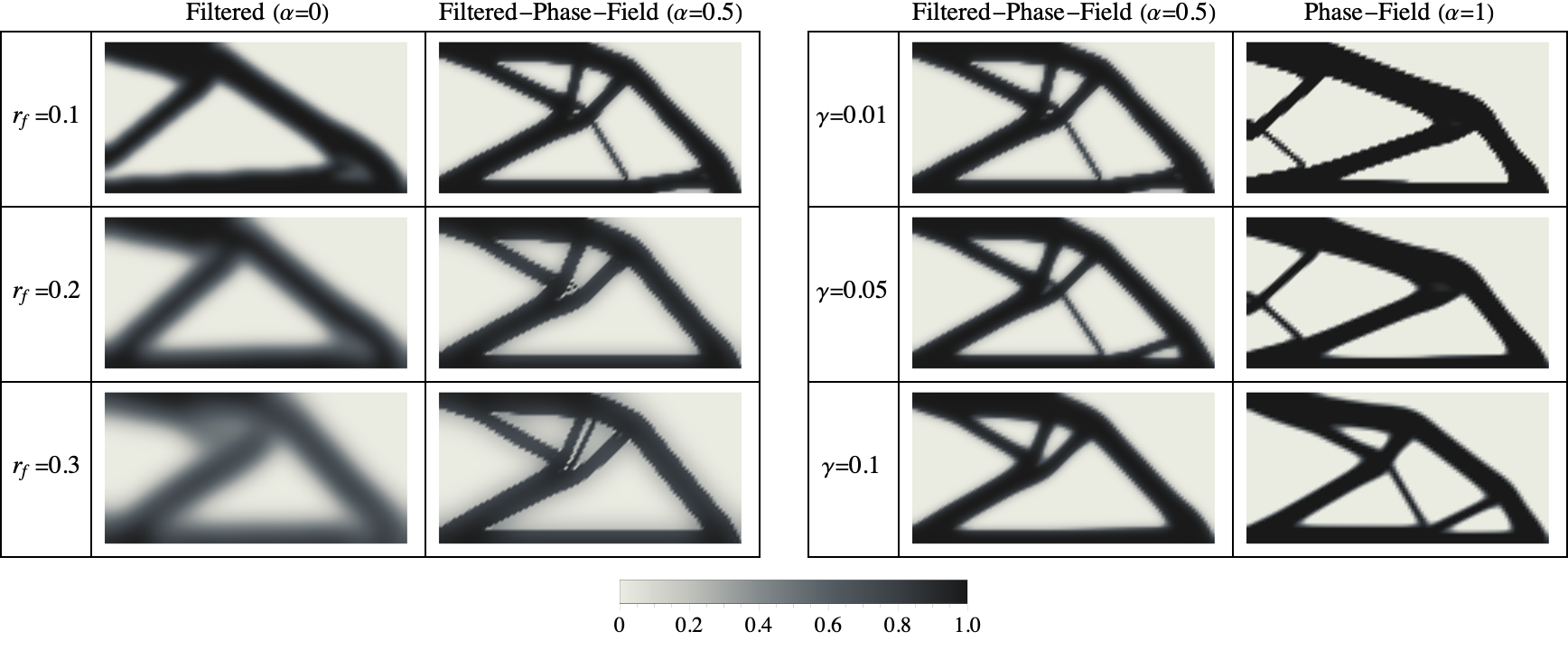}
\justifying

\caption{Cantilever beam topology optimization: optimal field $(\alpha \phi + \beta K \phi)$ obtained for $\alpha=1-\beta=0$ (F method), $\alpha=\beta=0.5$ (combined F/PF) and $\alpha=1-\beta=1$ (PF method) employing different values of the filter radius $r_f$ (for F and F/PF, left) and of the phase-field parameter $\gamma$ (for F/PF and PF, right). Simulation parameters (if not
  differently specified): $ \eta =1$ N/m, $r_f=0.1$ m, and $\gamma=0.01$ m.
} \label{fig}
\end{figure}

The purpose of this note is \UUU to present \EEE a combination of
F and PF in a single method, \UUU which we term {\it combined} F/PF
  {\it method} in the following. \EEE Such a combination is meant to set ground to a flexible
approach, where the features of the two methods are blended. 
In particular, given the parameters $\alpha, \, \beta,
\, \gamma \geq 0$, we investigate the bilevel minimization problem
\begin{equation}
  \label{eq:0}
  \boxed{\min_{\phi \in \Phi} \left\{C(\phi,u)+\alpha P_\gamma (\phi): \ u =
    \argmin_{v \in U} \big(E(\alpha\phi + \beta K \phi,v)-C(\phi,v)\big)\right\}}
\end{equation}
where the PF functional $P_\gamma$ for $\gamma >0$ or $\gamma =0$ is given by
\begin{align*}
  P_\gamma (\phi) &=
                    \left\{
                    \begin{array}{ll}\disp
  \eta \left[ \frac{\gamma}{2}  \int_\Omega |\nabla \phi|^2 \, \d
  x + \frac{1}{\gamma} \int_\Omega \phi^2(1-\phi)^2 \, \d x\right]
                      &\quad \text{if}  \ \phi \in H^1(\Omega)\\[1mm]
                      \infty&\quad \text{otherwise},
                    \end{array}
                              \right.\\
   P_0 (\phi) &=
                    \left\{
                    \begin{array}{ll}
    \disp \frac{\eta \sqrt{2}}{3} \,   {\rm  Per}(\{\phi=1\})  
                      &\quad \text{if}  \ \phi \in BV(\Omega,\{0,1\})\\[3mm]
                      \infty&\quad \text{otherwise},
                    \end{array}
                              \right.
\end{align*}
respectively. \UUU Different choices of the parameters $\alpha$ and
$\beta$ correspond to different combinations of F and PF, the
reference choice being $\beta=1-\alpha$, see also Figure~\ref{fig}. In
the following, we resort in keeping
these two parameters independent, for the sake of maximal generality. \EEE

In the definition of $P_0$, the symbol $ {\rm
  Per}(\{\phi=1\})  $ stands for the perimeter in $\Omega$ of the finite-perimeter
set $\{\phi=1\}$ and corresponds to the total variation $|{\rm
  D}\phi|(\Omega)$ of the Radon measure ${\rm
  D}\phi$.   The constant $\sqrt{2}/3$ in the definition of $P_0$ is explicitly computed from the specific form of the
 second term in $P_\gamma$ as $\sqrt{2}/3=2\int_0^1\sqrt{2 s^2 (1-s)^2} \, \d
 s $. 
The positive parameter $ \eta  $ is a weight factor required for unit consistency (and properly scaling the compliance and the PF functional), while $\gamma$ governs the thickness of the phase-field.

Throughout the paper, we work with a  linear
compact operator $K : L^\infty(\Omega) \to L^1(\Omega)$. A typical
example would be
\begin{equation}
  (K\phi)(x) = \int_\Omega k(x,y)\, \phi(y)\, \d y\label{eq:K}
\end{equation}
with $k\in W^{1,1}(\Omega^2)$. \UUU Although not strictly needed for the analysis,
 one may ask $k$ to be symmetric and $\int_\Omega k(x,y)\, \d x =1$ for
 all $y\in \Omega$. By additionally assuming that $\beta=1-\alpha$, this would
 entail that pure phases are conserved by the filter, namely $K1=1$
 and $K0=0$. Moreover, one would have that the volume constraint is
 also conserved, as $\int_\Omega K\phi\, \d x = \int_\Omega \phi \,
 \d x =v_0$ for all $\phi\in \Phi$.
 \EEE
 
The above notation and assumptions are considered in the following,
with no further explicit mention. The choice of the space
dimension $3$ is motivated by the application only. The arguments can
be easily recast in any dimension. 

The choice $\alpha=0$, $\beta=1$ returns the problem 
\begin{equation*}
  \min_{\phi \in \Phi} \left\{C(\phi,u) : \ u =
    \argmin_{v \in U} \big(E(  K \phi,v)-C(\phi,v)\big)\right\}
\end{equation*}
corresponding to the F method, where the
filtering is enacted by the compact operator $K$.
The purely filtered model \UUU $\alpha =0$ \EEE has been already investigated
numerically in \cite{Bruns} and theoretically in
\cite{Borrvall,Bourdin}, see also
\cite{Lazarov,Svanberg,Wadbro,Wang}. The effect of boundaries in
connection with the filtering is discussed in \cite{Clausen,Wallin}. The
filter radius as an additional design variable is treated in \cite{Amir}.

By choosing $\beta=0$, $\alpha=1$ we instead obtain
\begin{equation*}
  \min_{\phi \in \Phi} \left\{C(\phi,u)+ P_\gamma (\phi): \ u =
    \argmin_{v \in U} \big(E(\phi,v)-C(\phi,v)\big)\right\}
\end{equation*}
which is nothing but the PF method. 
The purely phase-field approach \UUU $\beta=0$ \EEE has been introduced in
\cite{Bourdin2,Burger} and developed in \cite{Blank,Garcke} for
multimaterials. Numerical investigations are in \cite{Dede,Takazawa}.
See also \cite{Auricchio} and for some extension to graded materials and
\cite{Almi,Almi2}  for some extension to graded materials and
elastoplasticity, respectively.


Both the F and the PF method are usually regarded as robust and
efficient. Still, their performance can significantly vary in specific
cases. Compared with its non-filtered version, the F method does not
show the occurrence of so-called {\it checkerboard modes}, namely fine
oscillations of solid and void,  \cite{Sigmund2007}. 
%
%
 However, the geometry of the optimal shape can be affected by a too coarse
filter, also leading to disconnections when defining a cut-off level
set for the definition of the solid domain. This issue is also
referred to as the occurrence of {\it grey transition regions}. 
In addition,  the geometry itself might depend on the filter radius
\UUU $r_f$, \EEE which for instance determines the size of the details
defined by the minimum feature sizes appearing in the final structure,
\cite{Sigmund2007}. \UUU Eventually, the choice of the threshold 
identifying the solid from the gray-scale computational output may strongly influence the
final topology. \EEE

\UUU The \EEE PF method does not require the introduction of
a filter, thus overcoming the  afore-mentioned 
issues. \UUU At the same time, it generally shows more robustness with
respect to the choice of the thresholding from the gray scale. On the
other hand, \EEE  the topology of the optimal shape is strongly
affected by the  choice of the  phase-field parameter $\gamma$, \UUU
which is usually just heuristic, see Figure \ref{fig}. \EEE 

The combined F/PF method \eqref{eq:0} seems in some cases to be able
to mitigate the criticalities of the underlying pure methods. \UUU In sharp
contrast with the pure F and PF methods, in the combined F/PF method
\UUU the optimal topology seems to be little
affected by the \UUU choice of the radius $r_f$  \EEE of the filter
and by  \UUU that of \EEE the phase-parameter $\gamma$. \UUU We find
this feature particularly important, for the parameters $r_f$ and
$\gamma$ are merely user-defined. The enhanced
robustness of the combined F/PF method with respect
to different choices of  $r_f$ and $\gamma$ is probably
the most interesting feature of this approach.

In order to illustrate the robustness of the combined F/PF method, we
present \UUU in  Figure  \ref{fig} a first \EEE numerical case study. We \EEE address a 2D cantilever beam: a rectangular
design region \UUU ($2 \times 1$ m, finite element size:
  $0.02 \times 0.02$ m) is clamped on the left side and   \UUU  loaded
  with a constant surface normal (downward) traction of 1 MPa acting
  upon the rightmost 10\% of the bottom side. \EEE   The topology optimization problem is solved by employing
a finite element space discretization in a Lagrangian formulation (see
the following Sections \ref{sec:space-discr} and
\ref{sec:lagr-form}).  Minimization is tackled via an 
Allen-Cahn gradient-based approach, adding a fixed global constraint
on the final volume (i.e., $\int_{\Omega}\phi=\bar{v}\,
 |\Omega|$   for some $\bar{v}\UUU = 0.4 \EEE$) and  suitably penalizing   the constraint
$\phi\in [0,1]$, see \cite{Marino}. The filter in  equation 
\eqref{eq:K} is built by means of {\it radial basis functions} with
support and observation points located in the middle of finite
elements (see, e.g., \cite{Shi}) and radius  $r_f$. Furthermore, the
functional dependence $\phi \mapsto \mathbb{C}(\phi)$ is introduced
with a classical SIMP power-law expression with power index equal to 3
and void $10^{-3}$ times softer than the solid material,
\cite{Sigmund2007}. \UUU More precisely, by letting the material constants of the
solid be $E=10$ GPa (Young's modulus) and $\nu=0.25$ (Poisson's
ratio), we let  $\mathbb{C}(\phi)=\mathbb{C}_0 +
\phi^3(\mathbb{C}_1-\mathbb{C}_0)$ with
$$    \mathbb{C}_1=\frac{E\nu}{1-\nu^2}\mathbb{I}_2\otimes
\mathbb{I}_2+ \frac{E}{1+\nu} \mathbb{I} \quad \text{and} \quad \mathbb{C}_0=10^{-3}\mathbb{C}_1$$
where $\mathbb{I}_2$ is the identity 2-tensor.  
By choosing $\beta=1-\alpha$, we report in  Figure  \ref{fig} \EEE the optimal shapes obtained for the pure F method ($\alpha=0$, $\beta=1$), the pure PF method ($\alpha=1$, $\beta=0$) and the combined F/PF method ($\alpha=\beta=0.5$) employing different values of the filter radius $r_f$ (for F and F/PF) and of the phase-field parameter $\gamma$ (for F/PF and PF).

\UUU Let us firstly compare the first two columns of  Figure
\ref{fig}, respectively corresponding to the F and the F/PF method, 
for different values of the filter radius $r_f$.  
The optimal phase distribution $\alpha \phi + \beta K \phi$
obtained with the F method depends on the radius $r_f$ of the
filter. Moreover, the choice  of the threshold level 
\EEE used to define the solid highly affects
the final geometry of the structure, possibly leading to
disconnections for coarse filters. \UUU 
The combined F/PF shows to be less
influenced by the filter radius $r_f$, as In fact, the main
geometrical features of the final topology are insensitive to the
choice of this parameter. Moreover, even employing filters
with large radii, grey transition regions are significantly less
present in the F/PF solution than in the corresponding F solution,
thus leading to a minor risk of creating disconnections when setting
the threshold \EEE for the definition of the solid. Remarkably, the
minimum feature sizes appearing in the final topology remain detailed
and small compared to the filter radius. It is noteworthy that no
post-processing technique has been employed for showing the solution
fields in Figure~\ref{fig}. Such techniques (see, e.g.,
\cite{Sigmund2007}) would allow to treat the issues of the F method,
 at the price of introducing  
an (undesirable) dependency  of  the final solution 
on user choices.

\UUU We now compare the last two columns of  Figure
\ref{fig}, respectively corresponding to the F/PF and the PF method, for varying values of the user-defined parameter $\gamma$. Optimal shapes from the PF method are highly affected
by the value of $\gamma$. On the contrary, final topologies from the
F/PF approach seem to be much less sensitive to the different choices
of the parameter. As $\gamma$ is usually fixed heuristically, we find the
robustness the F/PF method  with respect to this parameter 
particularly valuable. 
\EEE

In this note, we focus on the analytical aspects of method \eqref{eq:0}, while numerical and algorithmic considerations, as well as some
simulation campaigns, will be presented in a forthcoming publication.  At first, we prove that optimal shapes exist for any
choice of parameters $\alpha$, $\beta$, and $\gamma$  (Theorem
\ref{thm:ex}).   This existence proof is closely
reminiscent of that of \cite{Bourdin}.

We then prove that optimal shapes depend continuously (up to
subsequences) on the
parameters. In particular, letting $\phi_n$ be optimal for parameters
$(\alpha_n,\beta_n,\gamma_n)$, if $(\alpha_n,\beta_n,\gamma_n)\to
(\alpha,\beta,\gamma)$  then (up to not
relabeled subsequences) we have that $\phi_n \to \phi$, where $\phi$
is optimal for the parameters $(\alpha,\beta,\gamma)$ (Corollary
\ref{cor}). This convergence hinges on a more general variational
approximation result of $\Gamma$-convergence type (Theorem
\ref{thm:gamma}).

Problem \eqref{eq:0} is space-discretized by means of a Galerkin
method in Section~4. In particular, we discuss finite-dimensional
approximations of \eqref{eq:0} and prove the existence of
approximating optimal $\phi_h$. The convergence, up to subsequences,
of $\phi_h$ to solutions to \eqref{eq:0} is then recovered (Theorem \ref{thm:space}).

Eventually, we show in Theorem \ref{thm:lagrange} that \UUU the
stationary points of the \EEE bilevel minimization \UUU functional
$\phi \mapsto C(\phi,u) + \alpha P_\gamma(\phi)$ under the equilibrium
constraint of \eqref{eq:0} \EEE
can be equivalently tackled by finding stationary points
of the  Hu-Washizu-type {\it Lagrangian}
\begin{align*}
  L(\phi,u, e, \sigma) &= C(\phi,u)+\frac{\alpha}{2}P_\gamma(\phi) -
                         \frac12\int_\Omega \bC(\alpha \phi + \beta K
                         \phi) e : e \, \d x  \\
  &\quad + \int_\Omega \sigma: (e -
    \epsi(u))\, \d x.
    \end{align*}
 This equivalent formulation seem new in this context and allows for an
efficient numerical treatment of the topology optimization problem,  that is currently under investigations from the authors.

\section{Existence}

The focus of this section is on proving the existence of solutions to
problem \eqref{eq:0}. We have the following.

\begin{theorem}[Existence]\label{thm:ex}
  \UUU Problem \EEE \eqref{eq:0} admits a solution.
\end{theorem}

Before proving the result, let us collect notation on
the equilibrium problem \eqref{eq:a}-\eqref{eq:3}. Recall that for
any $  \phi \in \Phi$ one can find a unique $u \in U$ such that
  $$ E(\alpha\phi + \beta K \phi,u) - C(\phi,u) \leq  E(\alpha\phi +
  \beta K \phi,\haz u) - C(\phi,\haz u)\quad \forall \haz u \in U $$
  or, equivalently,
\begin{equation}\label{eq:sol} \int_\Omega \bC(\alpha \phi + \beta K
  \phi)\epsi(u):\epsi(v) \, \d x= \int_\Omega \phi f \cdot v \,
  \d x + \int_{\Gamma_{\rm N}} g \cdot f \, \d \Gamma \quad \forall v
  \in U.
  \end{equation}
  This allows to define a {\it solution operator} $S_{\alpha\beta} : \Phi \to U$ as
  $S_{\alpha\beta}(\phi)=u$.

  Standard estimates for the linear elastic system and the
  nondegeneracy of $\bC$ entail that
  $$c_0 \| \epsi (u)\|_{L^2}^2 \leq \| f \|_{L^2(\Omega;\Rz^3)} \| u
  \|_{L^2(\Omega;\Rz^3)} + \| g \|_{L^2( \Gamma_{\rm N};\Rz^3)} \| u
  \|_{L^2( \Gamma_{\rm N};\Rz^3)}.\label{eq:b0}$$
  Hence, an application of the Korn inequality ensures that $u$ 
  is bounded in $H^1(\Omega;\Rz^3)$. In particular,
  \begin{equation}
    \|S_{\alpha\beta}(\phi)\|_{H^1(\Omega;\Rz^3)} \leq C, \label{eq:b}
 \end{equation}
 where
   the constant $C$ depends on $\Omega$, $c_0$, $\| f
   \|_{L^2(\Omega;\Rz^3)} $, and $ \| g \|_{L^2( \Gamma_{\rm
       N};\Rz^3)}$ but is independent of $\phi$. This in turn implies
   that, the
   compliance term $C(\phi,S_{\alpha\beta}(\phi))$ is bounded,
   independently of $\phi\in \Phi$.

   By using the solution operator $S_{\alpha\beta}$ one can equivalently reformulate problem
  \eqref{eq:0} as
  \begin{equation*}
    \min_{\phi\in\Phi} \big\{C(\phi,S_{\alpha\beta}(\phi)) + \alpha P_\gamma(\phi) \big\}.
  \end{equation*}

\begin{proof}[Proof of Theorem \ref{thm:ex}]  
  Let $\phi_n$ be a minimizing sequence for $\phi \mapsto
  C(\phi,S_{\alpha\beta}(\phi)) + \alpha P_\gamma(\phi)$ and let $u_n =
  S_{\alpha\beta}(\phi_n)$. As $\phi_n\in [0,1]$ almost everywhere and $u_n$ is bounded in $H^1(\Omega;\Rz^3)$ by \eqref{eq:b},  by passing to some not
  relabeled subsequence we have that $\phi_n \wstarto \phi$ in
  $L^\infty(\Omega)$ and $u_n \weak u $ in $H^1(\Omega;\Rz^3)$. By
  compactness we also have that  $u_n \to u $ in
  $L^{6-}(\Omega;\Rz^3)$. In particular, we have that $f\cdot u_n \to
  f\cdot u$ in $L^{3/2-}(\Omega)$ and we can conclude that
  \begin{equation}\label{eq:L0}
C(\phi_n,u_n) \to C(\phi,u).
\end{equation}
\UUU We now proceed by considering separately the cases $\alpha >0$ and
$\alpha=0$.

{\it Case $\alpha >0$:} The \EEE boundedness of
  $ C(\phi_n, S_{\alpha\beta}(\phi_n)) + \alpha P_\gamma(\phi_n) $  and the fact that
  $C(\phi_n,S_{\alpha\beta}(\phi_n))$ is bounded  entail that $\phi_n$ is bounded in $H^1(\Omega)$ (if
  $\gamma>0$) or $BV(\Omega)$ (if $\gamma=0$). In both cases, $\phi_n$
  is precompact in $L^1(\Omega)$ so that one can assume with no
  loss of generality (or extract again, without relabeling) that $\phi_n \to \phi$ in
  $L^1(\Omega)$ as well. Hence, $\alpha \phi_n + \beta K \phi_n \to \alpha \phi + \beta K \phi  $ in
  $L^1(\Omega)$,
  \UUU independently of the value of $\beta \geq 0$. Due \EEE to the Lipschitz continuity and the
  boundedness of $\bC$ we hence have  that 
\begin{equation}\label{eq:L1}
\bC(\alpha \phi_n +
  \beta K \phi_n) \to \bC(\alpha \phi + \beta K \phi ) \quad \text{in}
  \ 
  L^q(\Omega;\Rz^{3\times 3\times 3\times 3 }), \ \forall q \in
  [1,\infty).
\end{equation}
 Since $\bC$ is uniformly
  positive, we equivalently have that $\bC(\alpha \phi_n +
  \beta K \phi_n)^{1/2} \to \bC(\alpha \phi + \beta K \phi )^{1/2}$ in
  $L^q(\Omega;\Rz^{3\times 3\times 3\times 3 })$ for all $q \in
  [1,\infty)$, where the superscript ${1/2}$ denotes the square-root tensor. This
  convergence, as well as the boundedness of $\bC$, allows us to
  conclude that 
\begin{equation}\label{eq:L2}\bC(\alpha \phi_n +
  \beta K \phi_n)^{1/2}\epsi(u_n)  \weak \bC(\alpha \phi + \beta K \phi
  )^{1/2} \epsi (u)\quad \text{in} \ 
  L^2(\Omega;\Rz^{3\times 3 }).
\end{equation}
 Let now $ u^* \in U$ and choose $u_m^*\in U  \cap
 W^{1,\infty}(\Omega,\Rz^3)$ with $ u_m^* \to u^*$ in
 $H^1(\Omega;\Rz^3)$. Making use of \eqref{eq:L1}-\eqref{eq:L2} we
 obtain
 \begin{align}
  & E(\alpha \phi + \beta K \phi,u) - C(\phi,u) \leq \liminf_{n\to
   \infty} \Big ( E(\alpha \phi_n + \beta K \phi_n,u_n) -
   C(\phi_n,u_n) \Big)\nonumber\\
&\quad \leq \liminf_{n\to
   \infty} \Big( E(\alpha \phi_n + \beta K \phi_n,u_m^*) -
   C(\phi_n, u_m^*) \Big)\nonumber\\[1mm]
&\quad = E(\alpha \phi + \beta K \phi, u_m^*) - C(\phi, u_m^*). \label{arg}
 \end{align}
Passing to the limit for $m\to \infty$, we have proved that $u = S_{\alpha\beta}(\phi)$.

\UUU Recall \EEE that we
have that $\phi_n \to \phi$ in
  $L^1(\Omega)$.
Moreover, since $\phi_n\in [0,1]$ a.e., this
  entails that $\phi_n \to \phi$ in
  $L^q(\Omega)$ for all $q \in [1,\infty)$.
In particular, if $\gamma>0$ one can pass the term $\int_\Omega\phi_n^2(1-\phi_n)^2
  \, \d x $ to the limit and check that 
$P_\gamma(\phi) \leq \liminf P_\gamma(\phi_n)$. If $\gamma=0$ one uses
the lower semicontinuity if the perimeter in $BV$ to get $P_0(\phi)
\leq \liminf P_0(\phi_n)$. As we have already checked that
$C(\phi_n,u_n) \to C(\phi,u)$ we conclude that 
$$C(\phi,u) + \alpha P_\gamma(\phi)  \leq \liminf_{n\to \infty} \Big(C(\phi_n,u_n) +\alpha
P_\gamma(\phi_n)\Big)$$
for all $\gamma\geq 0$ and $\phi$ is a solution of
problem \eqref{eq:0}.

\UUU 
{\it Case $\alpha =0$:} As $\phi_n \wstarto \phi$ in
$L^\infty(\Omega)$ and $K:L^\infty(\Omega) \to L^1(\Omega)$ is compact, we have that
$\beta K\phi_n \to \beta K \phi$ in $L^1(\Omega)$, independently of the
value of $\beta\geq 0$. Arguing as in \eqref{eq:L1} we have
\begin{equation}\label{eq:L10}
\bC(\beta K \phi_n)^{1/2} \to \bC(  \beta K \phi )^{1/2}\ \ \text{in}
\ \ L^q(\Omega;\Rz^{3\times 3\times 3\times 3 }), \ \forall q \in
  [1,\infty).
\end{equation}
This entails that 
\begin{equation}\label{eq:L2}\bC( 
  \beta K \phi_n)^{1/2}\epsi(u_n)  \weak \bC( \beta K \phi
  )^{1/2} \epsi (u)\quad \text{in} \ 
  L^2(\Omega;\Rz^{3\times 3 }).
\end{equation}
Repeating the argument in \eqref{arg} for $\alpha =0$ we get 
 \begin{align*}
  & E( \beta K \phi,u) - C(\phi,u) \leq E(\beta K \phi, u^*) - C(\phi,
    u^*)\quad \forall u^* \in U
 \end{align*}
 so that $u = S_{0\beta}(\phi)$. As
 \eqref{eq:L0} holds, this concludes the proof that  $\phi$ is a solution of
problem \eqref{eq:0}. \EEE
\end{proof}

\section{$\Gamma$-convergence and parameter asymptotics}

In this section, we provide asymptotic results  in relation
with limits in the parameters. We argue within the classical frame of
$\Gamma$-convergence \cite{DalMaso}. To this aim, it is notationally advantageous to
incorporate constraints into the definitions of the functionals by letting $G_{\alpha\beta\gamma}: L^\infty(\Omega) \times U \to
(-\infty,\infty]$ be defined as 
$$G_{\alpha\beta\gamma}(\phi,u) =
\left\{
  \begin{array}{ll}
    C(\phi,u) + \alpha P_\gamma(\phi)&\quad \text{if} \ \phi \in \Phi
      \ 
    \text{and} \ u = S_{\alpha\beta}(\phi)\\
\infty&\quad \text{otherwise}.
  \end{array}
\right.
 $$ 
Our $\Gamma$-convergence result reads as follows.
\begin{theorem}[$\Gamma$-convergence]\label{thm:gamma}
Let $(\alpha_n,\beta_n,\gamma_n) \to
  (\alpha,\beta,\gamma)$ with $\alpha >0 $ or $ \beta >0$. Then,
  $G_{\alpha_n\beta_n\gamma_n} \to G_{\alpha\beta\gamma}$ in the sense
  of $\Gamma$-convergence with respect to the weak$\ast$ topology of
  $ L^\infty(\Omega)\times U$.
\end{theorem}

\begin{proof} In order to prove the $\Gamma$-convergence
  $G_{\alpha_n\beta_n\gamma_n} \to G_{\alpha\beta\gamma}$  we check
  below the corresponding $\liminf$-inequality and we exhibit
  recovery sequences \cite{DalMaso}.
 
  {\it Liminf-inequality.} Let $(\phi_n,u_n)$ be given in such a way
  that $(\phi_n,u_n)\wstarto (\phi,u)$ in $L^\infty(\Omega)\times U$
  and assume with no loss of generality that 
\begin{equation}\label{eq:G}
\sup_{n\in \Nz}
  G_{\alpha_n\beta_n\gamma_n} (\phi_n,u_n)<\infty.
\end{equation}
 In particular,
  $\phi_n \in \Phi$ and $u_n =S_{\alpha_n\beta_n}(\phi_n)$ for all
  $n$, and, by possibly extracting without relabeling, we can assume that
  $u_n \to u$ in $L^{6-}(\Omega;\Rz^3)$. 

\UUU Let us now proceed by distinguishing cases.

{\it Case $\alpha>0$ and $\gamma>0$:} We have that $\phi_n$ are weakly precompact in $H^1(\Omega)$, hence strongly precompact in
$L^1(\Omega)$. This entails that 
\begin{equation}\label{eq:nun}
\alpha_n \phi_n + \beta_n K\phi_n \to \alpha
\phi + \beta K\phi \quad \text{in} \  L^1(\Omega)
\end{equation}
along a not relabeled subsequence, 
independently of the value of $\beta\geq 0$. Moving from \eqref{eq:nun},   the proof of
Theorem \ref{thm:ex} can be replicated verbatim, concluding that $u=
S_{\alpha\beta}(\phi)$.

In order to establish the $\Gamma$-$\liminf$ inequality
$G_{\alpha\beta\gamma } (\phi,u) \leq \liminf
G_{\alpha_n\beta_n\gamma_n } (\phi_n,u_n) $ we
just need to show that 
\begin{equation}
  \label{eq:liminf}
  C(\phi,u) + \alpha P_{\gamma}(\phi) \leq \liminf_{n \to \infty}
\Big( C(\phi_n,u_n) + \alpha_nP_{\gamma_n}(\phi_n)\Big).
\end{equation}
By following again the proof of Theorem \ref{thm:ex}, we obtain that
$C(\phi_n,u_n) \to C(\phi,u)$. 
One can  assume with
no loss of generality that $\alpha_n \geq \alpha/2$ and $\gamma_n \geq
\gamma/2>0$ for all $n$, so
that, possibly extracting without
relabeling, bound \eqref{eq:G} entails that $\phi_n \weak \phi$ in
$H^1(\Omega)$. Hence, $P_\gamma(\phi) \leq \liminf
P_{\gamma_n}(\phi_n)$ and we have that \eqref{eq:liminf} holds.

{\it Case $\alpha>0$ and $\gamma=0$:} We proceed as before in order to
check that  $u=
S_{\alpha\beta}(\phi)$. Inequality \eqref{eq:liminf} follows now by 
the classical Modica-Mortola result \cite{Modica}, yielding
$P_0(\phi) \leq \liminf P_{\gamma_n}(\phi_n)$.

{\it Case $\alpha=0$:} the compactness of $K$
suffices to conclude for 
\begin{equation}\label{eq:nun2}
\alpha_n \phi_n + \beta_n K\phi_n \to   \beta K\phi \quad \text{in} \  L^1(\Omega)
\end{equation}
independently of the value of $\beta\geq 0$. This again ensures that  $u=
S_{0\beta}(\phi)$.

The argument of  Theorem \ref{thm:ex} entails that 
$C(\phi_n,u_n) \to C(\phi,u)$ and the \linebreak
$\Gamma$-$\liminf$ inequality
\begin{align*}
&G_{0\beta\gamma } (\phi,u)  =  C(\phi,u) = \lim_{n\to \infty}
C(\phi_n,u_n) \\
&\quad \leq \liminf_{n\to \infty} \Big( C(\phi_n,u_n) + \alpha_n
P_{\gamma_n}(\phi_n) \Big)= \liminf_{n\to \infty}
G_{\alpha_n\beta_n\gamma_n } (\phi_n,u_n) 
\end{align*}
trivially follows, independently of the values of $\beta\geq 0$ and $\gamma\geq0$.
\EEE

{\it Recovery sequence.} Let $( \phi^*, u^*)\in L^\infty(\Omega)\times U$ with
$G_{\alpha\beta\gamma}( \phi^*, u^*)<\infty$ be given. We aim at
finding a recovery sequence $( \phi_n^*, u_n^*) \wstarto (
\phi^*, u^*)$ in $L^\infty(\Omega) \times U$ (at least) such that $G_{\alpha_n \beta_n \gamma_n}( \phi_n^*,  u
_n^*) \to G_{\alpha\beta\gamma}( \phi^*, u^*)$.

We distinguish here the two cases: $\gamma>0$
and $\gamma=0$.

\UUU {\it Case $\gamma>0$:} Independently of the values $\alpha, \,
\beta \geq 0$, we \EEE simply exploit pointwise convergence
and define $ \phi_n^* =  \phi^*$ and
$ u_n^* = S_{\alpha_n\beta_n}( \phi^*)$. As $\alpha_n  \phi^* +
\beta_n K  \phi^* \to \alpha \phi^* +
\beta K  \phi^* $ in $L^1(\Omega)$, it is  standard
to check that $ u_n^* \weak  u^*$ in $H^1(\Omega;\Rz^3)$,
so that one has 
$$C( \phi_n^* ,  u_n^*) + \alpha_n P_{\gamma_n}
( \phi_n^*) = C( \phi^* ,  u_n^*) + \alpha_n P_{\gamma_n}
( \phi^*) \to  C( \phi^* ,  u^*) + \alpha F_{\gamma}
( \phi^*)$$
and the convergence $G_{\alpha_n \beta_n \gamma_n}( \phi_n^*,  u
_n^*) \to G_{\alpha\beta\gamma}( \phi^*, u^*)$ follows.

\UUU {\it Case $\gamma=0$:} For all  $\alpha, \,
\beta \geq 0$, we \EEE resort again to the classical Modica-Mortola
construction \cite{Modica} in order to find $ \phi_n^* \to  \phi^*$
such that $P_{\gamma_n}( \phi_n^*) \to P_0y(
\phi^*)$. Correspondingly, we define again $ u_n^* =
S_{\alpha_n\beta_n}(\phi_n^*)$. As one again has that $\alpha_n  \phi_n^* +
\beta_n K  \phi_n^* \to \alpha \phi^* +
\beta K  \phi^* $ in $L^1(\Omega)$, one can still conclude that $
u_n^* \weak  u^*$ in $H^1(\Omega;\Rz^3)$. Hence, convergence $$C( \phi_n^* ,  u_n^*) + \alpha_n P_{\gamma_n}
( \phi_n^*)  \to  C( \phi^* ,  u^*) + \alpha F_{0}
( \phi^*)$$ holds
and $G_{\alpha_n \beta_n \gamma_n}( \phi_n^*,  u_n^*) \to G_{\alpha\beta0}( \phi^*, u^*)$ follows.
\end{proof}

The bound \eqref{eq:b} readily entails that the  functionals $G_{\alpha\beta\gamma}$ are equicoercive with respect
to the weak$\ast$ topology in
$L^\infty(\Omega)\times U$. The following is hence a straightforward
consequence of Theorem \ref{thm:gamma}.

\begin{corollary}[Parameter asymptotics]\label{cor}
  Let $\phi_{\alpha_n\beta_n\gamma_n}$ solve problem \eqref{eq:0} with
  parameters $(\alpha_n,\beta_n,\gamma_n)$. Moreover, let $(\alpha_n,\beta_n,\gamma_n) \to
  (\alpha,\beta,\gamma)$. Then, $\phi_{\alpha_n\beta_n\gamma_n} \wstarto
  \phi_{\alpha\beta\gamma}$ in $L^\infty(\Omega) $ up to a not relabeled
  subsequence, where 
$\phi_{\alpha\beta\gamma}$ solves  \eqref{eq:0} with
  parameters $(\alpha,\beta,\gamma)$. If $\alpha>0$ then
  $\phi_{\alpha_n\beta_n\gamma_n}$ converges also strongly in $L^q(\Omega) $ for all $q \in
  [1,\infty)$. \UUU Moroever, it converges \EEE weakly in $H^1(\Omega)$ if $\gamma>0$, and weakly$\ast$
  in $BV(\Omega)$ if $\gamma=0$.
\end{corollary}

\section{Space discretization} \label{sec:space-discr}

We describe now a space-discretization procedure via a
Galerkin method. Although
our approach is abstract, assumptions are modeled on the case of
conformal finite elements. Let $\Phi_h \subset \Phi$ and $U_h \subset U$ be two
families of finite-dimensional subspaces with $\Phi_h \subset \Phi_{h'}$
and $U_h \subset U_{h'}$ for $h'\leq h$, $\cup_h\Phi_h$ dense in
$\Phi$, and $\cup_hU_h$ dense in
$U$. We also assume that $\cup_hU_h$ is dense in
$W^{1,\infty}(\Omega;\Rz^3)$ with respect to its topology. These
assumptions are fulfilled by choosing $\Phi_h$ and $U_h$ as spaces
of piecewise polynomials of degree $0$ and $1$ on a given regular
triangulation of $\Omega$ (assume it to be a polygon) with mesh size
$h>0$.

For the sake of simplicity, we assume to be able to evaluate the
functionals $E$, $P_\gamma$, and $C$ exactly on $\Phi_h$ and
$U_h$. Note however that the analysis can be extended to the case of
approximating $E_h$, $F_{\gamma h}$, and $C_h$ at the expense of some
additional
notational intricacy only, see \cite[Sec. IV.27.4.2]{Ern}, as well as the classical
references \cite{Ciarlet,Strang}. On the contrary, we assume to be given a
family of linear and continuous operators $K_h: \Phi_h \to \Phi_h$
fulfilling the {\it continuous convergence} requirement
\begin{equation}
  \label{eq:cc}
  \phi_h \wstarto \phi \ \ \text{in} \ L^\infty(\Omega) \quad
  \Rightarrow \quad  K_h\phi_h \to K\phi \ \ \text{in} \ L^1(\Omega).
\end{equation}
The latter can be readily met in practice, if $K$ is chosen to have form
\eqref{eq:K}.

The space-discrete version of problem \eqref{eq:0} reads as follows
\begin{align}
 \nonumber
 & \min_{\phi_h \in \Phi_h } \bigg\{C(\phi_h ,u_h )+\alpha P_\gamma
   (\phi_h ): \\
  &\quad  u_h  =
    \argmin_{v_h  \in U_h } \big(E(\alpha\phi_h  + \beta K_h  \phi_h ,v_h )-C(\phi_h ,v_h )\big)\bigg\} \label{eq:0d}
\end{align}
Note that the latter makes sense, for $v_h   \mapsto
E(\alpha\phi_h  + \beta K_h  \phi_h ,v_h )-C(\phi_h ,v_h )$ admits a
unique minimizer in $U_h$ for all $\phi_h\in \Phi_h$ due to the
Lax-Milgram Lemma. In particular, this defines the discrete solution
operator $S_{\alpha\beta h}: \Phi_h \to U_h$ as $u_h = S_{\alpha\beta h}(\phi_h)$, allowing to
equivalently rewrite problem \eqref{eq:0d} as follows
\begin{align}
 & \min_{\phi_h \in \Phi_h } \big\{C(\phi_h ,S_{\alpha\beta h}(\phi_h) )+\alpha P_\gamma
   (\phi_h ) \big\}. \label{eq:1d}
\end{align}
The main result of this section is the following.

\begin{theorem}[Space discretization]\label{thm:space} For all $h>0$ problem \eqref{eq:1d} admits a
  solution $\phi_h$. There exists a positive constant $C$ depending on
  $c_0$, $\Omega$, $\| f \|_{L^2(\Omega;\Rz^3)}$, and $\| g
  \|_{L^2(\Gamma_{\rm N};\Rz^3)}$ but independent of $h$, $\alpha$,
  and $\beta$ such that $
  \| u_h \|_{H^1(\Omega;\Rz^3)} \leq C$, for $u_h = S_{\alpha\beta h}(\phi_h)$. As
  $h \to 0$, one can find not relabeled subsequences  such that $\phi_h
  \wstarto \phi$ in $L^\infty(\Omega)$ and $u_h \weak u$ in
  $H^1(\Omega;\Rz^3)$, where $\phi$ solves the limiting problem
  \eqref{eq:0} and $u=S_{\alpha\beta}(\phi)$.
 If $\alpha>0$ then
  $\phi_h$ converges also strongly in $L^q(\Omega) $ for all $q \in
  [1,\infty)$, weakly in $H^1(\Omega)$ if $\gamma>0$, and weakly$\ast$
  in $BV(\Omega)$ if $\gamma=0$.
\end{theorem}

\begin{proof}
  The existence of space-discrete solutions $\phi_h$ follows by the
  same argument as in Theorem \ref{thm:ex}. The situation is even
  simpler here, for the finite dimensionality of the problem
  entails that, for $h$ fixed, a minimizing sequence $\phi_{hn}$ for
  \eqref{eq:1d} is strongly compact in $L^\infty(\Omega)$ (and, for $\alpha>0$, in
  $H^1(\Omega)$ or $BV(\Omega)$, depending on $\gamma$).

Let now $\phi_h$ solve \eqref{eq:1d} and $u_h = S_{\alpha\beta h}(\phi_h)$. 
The bound on $\| u_h \|_{H^1(\Omega;\Rz^3)}$ can be obtained as in
\eqref{eq:b}. As  $\phi_h\in [0,1]$ a.e., as $h \to 0$ one can extract
(without relabelling) and have $\phi_h
  \wstarto \phi$ in $L^\infty(\Omega)$ and $u_h \weak u$ in
  $H^1(\Omega;\Rz^3)$. Moreover, the above bounds ensure that
  $C(\phi_h,u_h)$ is bounded, independently of $h$. Hence,
  if $\alpha>0$ one can assume that 
  $\phi_h$ converges also strongly in $L^q(\Omega) $ for all $q \in
  [1,\infty)$, weakly in $H^1(\Omega)$ if $\gamma>0$, and weakly$\ast$
  in $BV(\Omega)$ if $\gamma=0$.

In all cases, by using the continuous-convergence assumption
\eqref{eq:cc} we have that $\alpha \phi_h + \beta K_h \phi_h \to
\alpha \phi  + \beta K  \phi $ in $L^1(\Omega)$. By following
the argument of Theorem \ref{thm:ex}, we again obtain that 
\begin{align}
  &\bC(\alpha \phi_h + \beta K_h \phi_h ) \to \bC(\alpha \phi  + \beta
    K \phi  ) \quad \text{in} \
L^q(\Omega;\Rz^{3\times 3\times 3\times 3})\quad \forall q \in [1,\infty),\label{eq:conv1}\\
 &\bC(\alpha \phi_h + \beta K_h \phi_h )^{1/2}\epsi (u_h) \weak
\bC(\alpha \phi + \beta K \phi)^{1/2}\epsi (u)\quad \text{in} \
L^2(\Omega;\Rz^{3\times 3}). \label{eq:conv2}
\end{align}
Let now $v_m \in W^{1,\infty}(\Omega;\Rz^3)$ be given and approximate
it via $v_h\in U_h$ such that $v_h \to v_m$ in
$W^{1,\infty}(\Omega;\Rz^3)$ as $h\to 0$. Using the convergences \eqref{eq:conv1}-\eqref{eq:conv2} and the fact that $C(\phi_h, u_h) \to C(\phi , u)$ we deduce that
\begin{align*}
  &E(\alpha \phi + \beta K \phi,u) - C(\phi, u) \leq \liminf_{h\to 0}
    \Big( E(\alpha \phi_h + \beta K_h \phi_h,u_h) - C(\phi_h,
    u_h)\Big) \\
  &\quad \leq \liminf_{h\to 0}
    \Big( E(\alpha \phi_h + \beta K_h \phi_h,v_h) - C(\phi_h,
    v_h) \Big)= E(\alpha \phi + \beta K \phi,v_m) - C(\phi, v_m).
\end{align*}
Eventually, for all $v \in U$ one can find a sequence $v_m \in
W^{1,\infty}(\Omega;\Rz^3)$ such that $v_m \to v$ in
$H^1(\Omega;\Rz^3)$ as $m\to \infty$. By passing to the limit as $m\to
\infty$ in the above inequality ensures that $u \in S_{\alpha\beta}(\phi)$.

The last step of the proof consists in remarking that
$$C(\phi,u) + \alpha P_\gamma(\phi) \leq \liminf_{h\to 0} \Big( C(\phi_h,u_h) +
\alpha P_\gamma(\phi_h) \Big)$$
\UUU independently of the values of the parameters  $(\alpha,  \beta, \gamma)$. In
fact, one has that $\UUU \alpha  \EEE P_\gamma(\phi) \leq \liminf_{h\to 0}  
\alpha P_\gamma(\phi_h) $, for any $\alpha\geq 0$. \EEE Eventually, we have
checked that $\phi$ solves \eqref{eq:1d}.
\end{proof}

Before closing this section, let us remark that the above analysis can
be extended to include parameter asymptotics, in the same spirit of
Section 3.

\section{Lagrangian formulation} \label{sec:lagr-form}

The actual implementation of the bilevel minimization of problem
\eqref{eq:0} is computationally demanding. \UUU On the contrary, stationary
points of the bilevel minimization functional $\phi\mapsto
C(\phi,S_{\alpha\beta}(\phi)) + \alpha P_\gamma(\phi)$ can be efficiently tackled by
equivalently reformulating the \EEE problem in terms of stationarity of
the {\it Lagrangian} $L : \Phi\times U \times L^2(\Omega;\Rz^{3\times 3})
\times L^2(\Omega;\Rz^{3\times 3}) \to \Rz$ given by 
\begin{align*}
  L(\phi,u, e, \sigma) &= C(\phi,u)+\frac{\alpha}{2}P_\gamma(\phi) -
                         \frac12\int_\Omega \bC(\alpha \phi + \beta K
                         \phi) e : e  \, \d x  \\
  &\quad + \int_\Omega \sigma: (e -
    \epsi(u))\, \d x.
    \end{align*}

The fact that \UUU stationary points of $\phi\mapsto
C(\phi,S_{\alpha\beta}(\phi)) + \alpha P_\gamma(\phi)$ \EEE and (the first component of)
stationary points of the Lagrangian $L$ coincide was already used without
proof  in \cite{Marino} in the setting of the PF method ($\alpha =1$,
$\beta=0$). On the other hand, such {\it monolithic} formulation seems
to be new in the frame of the F method ($\alpha =0$,
$\beta=1$). 

Note that, for the purposes of simplifying the presentation, the
constraints $\phi\in [0,1]$ and $\int_\Omega \phi \, \d x =v_0$ are
neglected throughout this section. Our main result is the following.

\begin{theorem}[Lagrangian formulation]\label{thm:lagrange} \UUU
  $\phi$ is a stationary point of \EEE
  $C(\phi,S_{\alpha\beta}(\phi)) + \alpha P_\gamma(\phi)$ if and only if the Lagrangian $L$ is
  stationary at $(\phi,S_{\alpha\beta}(\phi), \epsi (S_{\alpha\beta}(\phi)), \bC(\alpha \phi + \beta K \phi)\epsi(S_{\alpha\beta}(\phi)))$.
\end{theorem}
 
\begin{proof} 
  By computing variations of \eqref{eq:sol} for $u=S_{\alpha\beta}(\phi)$ in direction 
 $\tilde \phi \in L^\infty(\Omega)$  (in the case with
  constraints we would require $\int_\Omega \tilde \phi \, \d x =0$ and $\phi
  + t \tilde \phi \in [0,1]$ a.e., for $t$ small enough), we get that
  \begin{align}
   & \int_\Omega \bC(\alpha \phi + \beta K
     \phi)\epsi(D S_{\alpha\beta}(\phi)[\tilde \phi]):\epsi(v) \, \d x \nonumber\\
    &\qquad + \int_\Omega D\bC (\alpha \phi + \beta K
  \phi) [\alpha \tilde \phi + \beta K
    \tilde \phi]\epsi(S_{\alpha\beta}(\phi)):\epsi(v) \, \d x\nonumber\\
    &\quad= \int_\Omega \tilde \phi f \cdot v \,
  \d x \quad \forall v
    \in U. \label{eq:var}
  \end{align}
Here,  $D S_{\alpha\beta}(\phi)[\tilde \phi]\in U$ and $D\bC (\alpha \phi + \beta K
  \phi) [\alpha \tilde \phi + \beta K
    \tilde \phi]$ are the Gateaux derivatives of $\phi\mapsto S_{\alpha\beta}(\phi)$
    and $\phi\mapsto  \bC(\alpha \phi + \beta K
     \phi) $ at
$\phi$ in direction $\tilde \phi$, respectively. 

Compute now the variations of $L$ at $(\phi,u,e,\sigma)$ in directions $(\tilde \phi, \tilde u , \tilde
e , \tilde \sigma)$ and get
\begin{align*}
  \delta_\phi L(\phi,u,e,\sigma)[\tilde \phi, \tilde u , \tilde
  e , \tilde \sigma]&=\int_\Omega \tilde \phi f \cdot u \, \d x
                      +\frac{\alpha}{2}\delta P_\gamma(\phi)[\tilde
                      \phi] \\
  &\quad -\frac12\int_\Omega D\bC(\alpha \phi +
                      \beta K \phi)[\alpha \tilde\phi + \beta K
                      \tilde\phi] e: e \, \d x,   \\
   \delta_u L(\phi,u,e,\sigma)[\tilde \phi, \tilde u , \tilde
  e , \tilde \sigma]&= \int_\Omega \phi f \cdot \tilde u \, \d x +
                      \int_{\Gamma_{\rm N}} g \cdot \tilde u \,\d
                      \Gamma + \int_\Omega \sigma : \epsi (\tilde u),
  \\
   \delta_e L(\phi,u,e,\sigma)[\tilde \phi, \tilde u , \tilde
  e , \tilde \sigma]&= - \int_\Omega \big(\bC (\alpha \phi +
                      \beta K \phi) e - \sigma\big): \tilde e\, \d x ,
  \\
   \delta_\sigma L(\phi,u,e,\sigma)[\tilde \phi, \tilde u , \tilde
  e , \tilde \sigma]&=   \int_\Omega \tilde \sigma: (e -
\epsi(u))\, \d x.   
\end{align*}
Note that  
\begin{align*}
  \delta_e L(\phi,u,e,\sigma) = 0 \quad & \Leftrightarrow \quad \sigma
                                          =  \bC (\alpha \phi +
                      \beta K \phi) e \quad \text{a.e.},\\
  \delta_\sigma L(\phi,u,e,\sigma) = 0 \quad & \Leftrightarrow \quad
                                                 e = \epsi (u) \quad \text{a.e}.
\end{align*}
Indeed, stationarity in $e$ and $\sigma$ deliver the constitutive
equation and the kinematic compatibility, respectively. On the other hand,
stationarity in $u$ corresponds to equilibrium. In particular,
  $u = S_{\alpha\beta}(\phi)$ turns out to be equivalent to $\delta_uL =0$,
$\delta_eL=0$, and $\delta_\sigma L=0$ at the point
$(\phi,u,e,\sigma)=(\phi,S_{\alpha\beta}(\phi), \epsi (S_{\alpha\beta}(\phi)), \bC(\alpha \phi +
\beta K \phi)\epsi(S_{\alpha\beta}(\phi)))$.

In order to conclude the proof, we compute the variation of $\phi
\mapsto C(\phi,S_{\alpha\beta}(\phi)) + \alpha P_\gamma(\phi)$ in direction $\tilde
\phi$ getting
\begin{align*}
 & \delta \big(C(\phi,S_{\alpha\beta}(\phi)) + \alpha P_\gamma(\phi) \big)[\tilde\phi] =
  \int_\Omega \tilde \phi f \cdot S_{\alpha\beta}(\phi)\, \d x + \int_\Omega  \phi
  f \cdot D S_{\alpha\beta}(\phi)[\tilde\phi]\, \d x \\
  &\quad + \int_{\Gamma_{\rm N}} 
  g \cdot D S_{\alpha\beta}(\phi)[\tilde\phi]\, \d \Gamma + \alpha \delta
  P_\gamma(\phi)[\tilde \phi].
\end{align*}
By using relation \eqref{eq:var} and setting  $u=S_{\alpha\beta}(\phi)$ we obtain that
\begin{align*}
  & \delta \big(C(\phi,S_{\alpha\beta}(\phi)) + \alpha P_\gamma(\phi) \big)[\tilde \phi]
  \\
  &\quad=
\int_\Omega \tilde \phi f \cdot S_{\alpha\beta}(\phi)\, \d x  + \int_\Omega \bC (\alpha \phi +
                      \beta K \phi) \epsi(DS_{\alpha\beta}(\phi)[\tilde \phi]) :
    \epsi(u)\, \d x + \alpha \delta
    P_\gamma(\phi)[\tilde \phi] \\
  &\quad =2\int_\Omega \tilde \phi f \cdot S_{\alpha\beta}(\phi)\, \d x  - \int_\Omega D\bC (\alpha \phi +
                      \beta K \phi)[\alpha \tilde \phi + \beta K
    \tilde \phi] \epsi(u) : \epsi(u) \, \d x \\
  &\qquad + \alpha \delta
    P_\gamma(\phi)[\tilde \phi]\\
  &\quad= 2 \delta_\phi L(\phi, u, e, \sigma )[\tilde \phi, \tilde  u,
    \tilde  e, \tilde \sigma ] 
\end{align*} 
and the assertion follows. 
\end{proof}


\section*{Acknowledgement}
F. Auricchio was partially supported by the Italian Minister of
University and Research through the project {\it A BRIDGE TO THE
  FUTURE: Computational methods, innovative applications, experimental
  validations of new materials and technologies} (No. 2017L7X3CS)
within the PRIN 2017 program abd by Regione Lombardia, regional law no. 9/2020, resolution no. 3776/2020.
M. Marino was partially supported by the Italian Ministry of
University and Research through the project {\it COMETA} within the
Program for Young Researchers {\it Rita Levi Montalcini} (year 2017)
and by Regione Lazio through the project {\it BIOPMEAT}
(No. A0375-2020-36756) within the framework {\it Progetti di Gruppi di Ricerca 2020} (POR FESR LAZIO 2014). 
I. Mazari is partially supported by the French ANR Project
ANR-18-CE40-0013-SHAPO on Shape Optimization and by the Project {\it Analysis and simulation of optimal shapes - application to life sciences} of the Paris City Hall.
U. Stefanelli is partially supported by 
the Austrian Science Fund (FWF) through projects F\,65, W\,1245,  I\,4354, I\,5149, and P\,32788, and by the OeAD-WTZ project CZ
01/2021. The authors have no relevant financial or non-financial interests to disclose.

\end{document}